\documentclass[12pt]{article}
\usepackage[margin=0.7in]{geometry}
\usepackage{amssymb,amsmath,amsthm,mathrsfs}
\usepackage{graphicx,epstopdf,color}

\allowdisplaybreaks
\usepackage[alphabetic]{amsrefs}
\usepackage [latin1]{inputenc}
\newcommand{\R}{\mathbb{R}}

\begin{document}

\title{Alessio Figalli's contributions to nonlocal minimal surfaces}

\author{Serena Dipierro\and Enrico Valdinoci}

\maketitle

{\scriptsize \begin{center} Department of Mathematics and Statistics\\
University of Western Australia\\ 35 Stirling Highway, WA6009 Crawley (Australia)\\
\end{center}
\bigskip

\begin{center}
E-mail addresses:
{\tt serena.dipierro@uwa.edu.au},
{\tt enrico.valdinoci@uwa.edu.au}
\end{center}
}
\bigskip\bigskip

Alessio has produced in his very intense career an extraordinary number of outstanding results
in an impressive variety of topics. Among the multifold research lines in which he
acted as a trailblazer, the one focused on nonlocal minimal surfaces
offered an excellent opportunity for Alessio to pioneer some of the first settlements
in a brand new subject of investigation and pave the way to a broad spectrum of future
research.\medskip

Nonlocal minimal surfaces are beautiful objects whose research combines
motivations and methods arising in different disciplines, including mathematical
analysis, differential geometry and mathematical physics, and a full understanding
of their complexity requires a truly cross-disciplinary and open-minded approach.
As it often occurs in the new research lines, to understand nonlocal minimal surfaces
one has to discover novel methodologies revealing the striking differences
with respect to the previous knowledge, get off the beaten path, and think differently.
Giving a full account of all the important progress that the field
of nonlocal minimal surfaces has recently experienced is a goal which goes well beyond
the purpose of this note, therefore we will simply focus here on some
of the very original and important contributions given by Alessio in this field.
Without any attempt of being exhaustive, other fundamental contributions
provided by other authors will be only tangentially discussed:
the reader who wants to dig more into the subject can consult the existing literature,
including a set of lecture notes coauthored by Alessio himself~\cite{MR3588123}.
\medskip

Nonlocal minimal surfaces have been introduced in~\cite{MR2675483}
as the outcome of a minimization problem involving a nonlocal notion of perimeter.
Roughly speaking, the energy functional takes into account the pointwise interactions
of a set with its complement. These interactions are weighted by a kernel which
is invariant under translations and rotations, and which is self-similar after scaling.
One can also take into account the contributions
of this energy functional with respect to a given reference domain~$\Omega\subseteq\R^n$
(that we take with smooth boundary for the sake of simplicity).
In this case the energy contributions that involve the interactions
of points lying in the complement of~$\Omega$ are omitted from the functional.

\begin{figure}
    \centering
    \includegraphics[width=8.9cm]{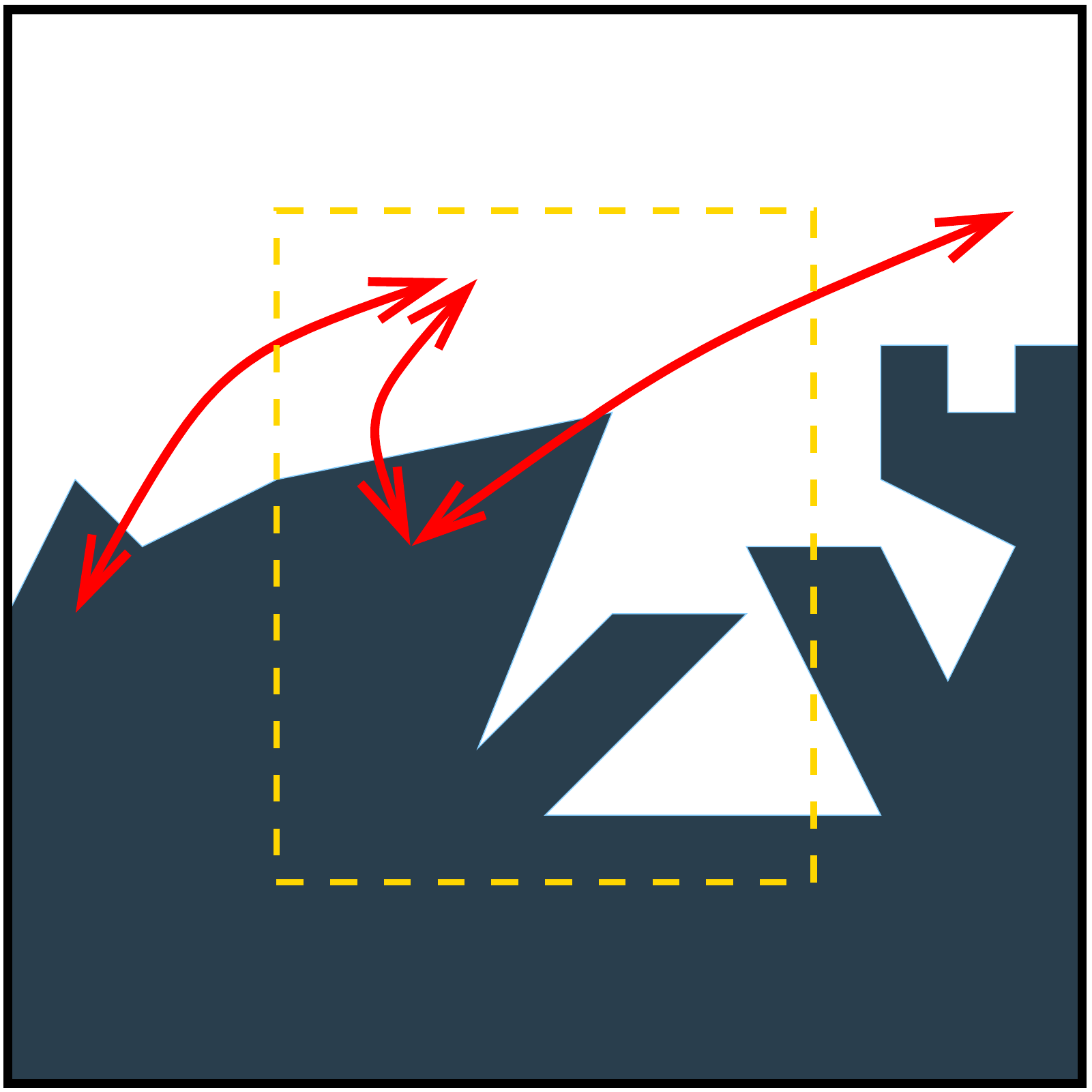}
    \caption{\em {{Pointwise interactions defining the $s$-perimeter.}}}
    \label{DOMA}
\end{figure}

More precisely, fixed~$s\in(0,1)$, one defines the $s$-perimeter of a (measurable)
set~$E\subseteq\R^n$ with respect to the reference domain~$\Omega$ as
\begin{equation}\label{INOM} {\rm Per}_s(E,\Omega):=
I_s(E\cap\Omega, E^c\cap\Omega)+I_s(E\cap\Omega, E^c\cap\Omega^c)+
I_s(E\cap\Omega^c, E^c\cap\Omega),\end{equation}
where we used the superscript ``$c$'' to denote complementary sets in~$\R^n$,
and~$I_s(\cdot,\cdot)$ represents the set interaction given by
\begin{equation}\label{Is} I_s(A,B):=s(1-s)\iint_{A\times B}\frac{dx\,dy}{|x-y|^{n+s}},\qquad{\mbox{ for all
disjoint subsets $A$ and $B$ of~$\R^n$,}}\end{equation}
see Figure~\ref{DOMA}.

In spite of its remarkable structural simplicity, codifying the most fundamental
geometric property of a given set, the $s$-perimeter turns out to be one
of the most difficult objects to fully understand. On the other hand,
it provides a large amount of information on several models of concrete
interest, including long-range phase transitions,
spin models and image reconstruction.\medskip

The factor $s(1-s)$ has been included in the interaction functional in~\eqref{Is}
for normalization purposes. In this way, one has that
if~$E$ has finite classical perimeter in a
neighborhood of~$\Omega$, then
$$\lim_{s\nearrow1}{\rm Per}_s(E,\Omega) = \omega_{n-1}\, {\rm Per}(E,\overline\Omega),$$
where~$\omega_n$ is the volume of the unit ball in~$\R^n$, and ``{\rm Per}''
denotes the classical perimeter. In this sense, up to a dimensional constant,
the $s$-perimeter recovers the classical perimeter as~$s\nearrow1$.

The limit of the $s$-perimeter as~$s\searrow0$ is somehow more complicated,
and it has been investigated in full details in~\cite{MR3007726}.
In this work there are explicit examples of smooth
sets for which such a limit does not even exist, and, in general, it is shown that
the existence of the limit of the $s$-perimeter as~$s\searrow0$
is strictly related to the existence of the following limit:
\begin{equation}\label{ZETA}
\zeta(E) := \lim_{s\searrow0}
\frac{s}{n\omega_n}\int_{E\setminus B_1}\frac{dx}{|x|^{n+s}}.
\end{equation}
Roughly speaking, the quantity~$\zeta(E)$ measures the ``mass''
of the set~$E$ at infinity, weighted by the kernel.
We observe that~$\zeta$ is monotone with respect to set inclusion.
Also, it is apparent that~$\zeta(\varnothing)=0$ and~$\zeta(\R^n)=1$.
Furthermore, $\zeta$ evaluated at a half-space is exactly~$1/2$
and, more generally, $\zeta$ evaluated at a cone gives the opening of the cone itself, that is
if~$E:=\{ tp,\;p\in\Sigma,\;t>0\}$ for some~$\Sigma\subseteq S^{n-1}$, then
\begin{eqnarray*}
\zeta(E) &=& \lim_{s\searrow0}
\frac{s}{{\mathcal{H}}^{n-1}(S^{n-1})}
\iint_{(\rho,\omega)\in
(1,+\infty)\times \Sigma}\frac{\rho^{n-1}}{\rho^{n+s}}\,d\rho\,d{\mathcal{H}}^{n-1}(\omega)\\
&=&\frac{{\mathcal{H}}^{n-1}(\Sigma)}{{\mathcal{H}}^{n-1}(S^{n-1})}.
\end{eqnarray*}
It turns out that
we can consider~$\zeta(E)$ as a ``convex parameter'',
and it is proved in~\cite{MR3007726} that
if~$E$ has finite~$s_0$-perimeter in~$\Omega$
for some~$s_0\in(0,1)$
and the limit in~\eqref{ZETA} exists, then
the limit as~$s\searrow0$
of the $s$-perimeter of~$E$ in~$\Omega$ also exists, and
\begin{equation}\label{AS0} \lim_{s\searrow0} {\rm Per}_s(E,\Omega)=(1-\zeta(E))\,|E\cap\Omega|+
\zeta(E)\,|E^c\cap\Omega|,
\end{equation}
where we denoted by~$|\cdot|$ the volume of a set.

In the special case in which~$|E\cap\Omega|=|E^c\cap\Omega|$,
it is also shown in~\cite{MR3007726} that~\eqref{AS0}
is true independently on the existence of the limit in~\eqref{ZETA},
since in this case
\begin{equation*} \lim_{s\searrow0} {\rm Per}_s(E,\Omega)=|E\cap\Omega|=|E^c\cap\Omega|.
\end{equation*}
On the other hand, if~$|E\cap\Omega|\ne|E^c\cap\Omega|$
the existence of the limit in~\eqref{AS0} is shown to be equivalent to that
in~\eqref{ZETA}.

\begin{figure}
    \centering
    \includegraphics[width=6.9cm]{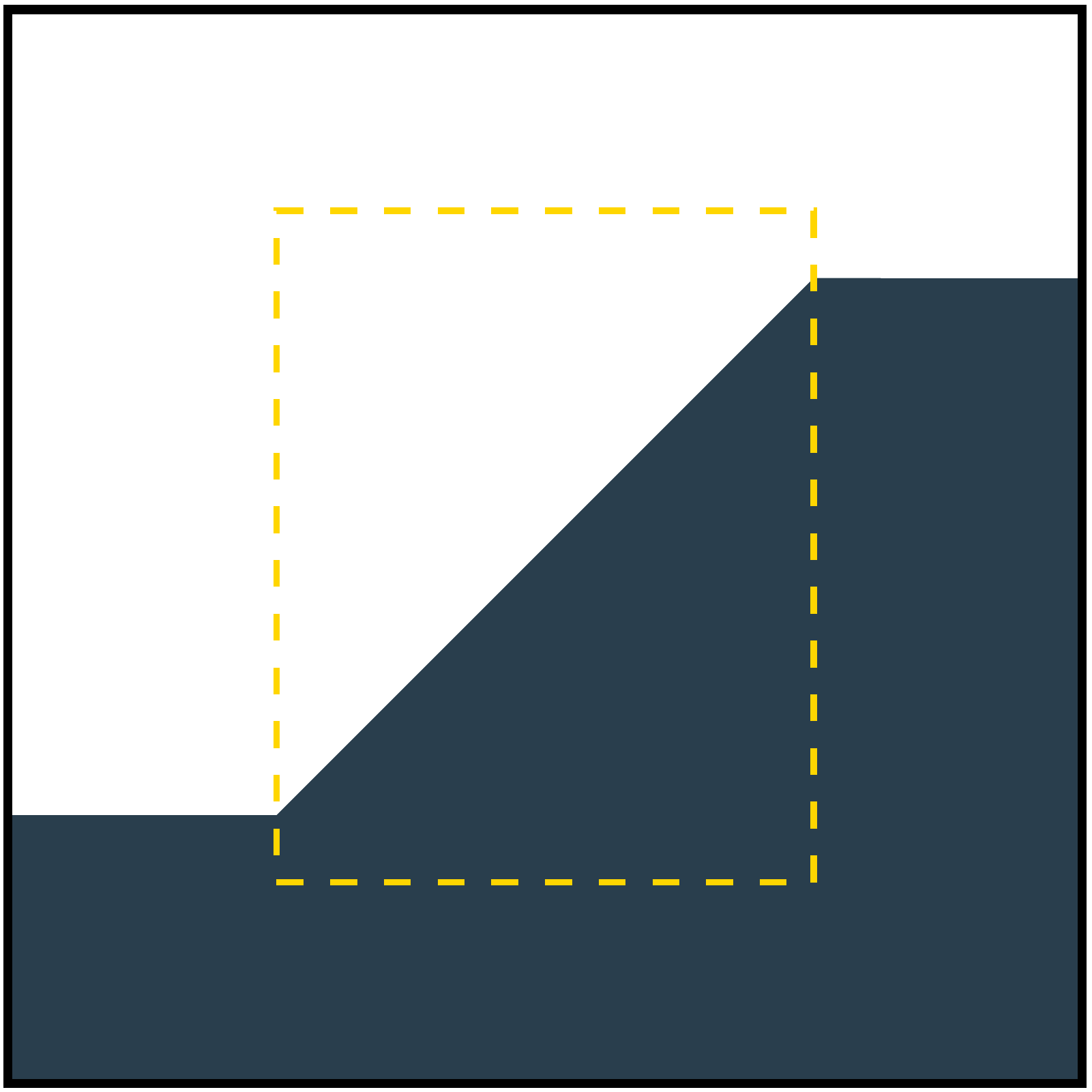} $\qquad$
    \includegraphics[width=6.9cm]{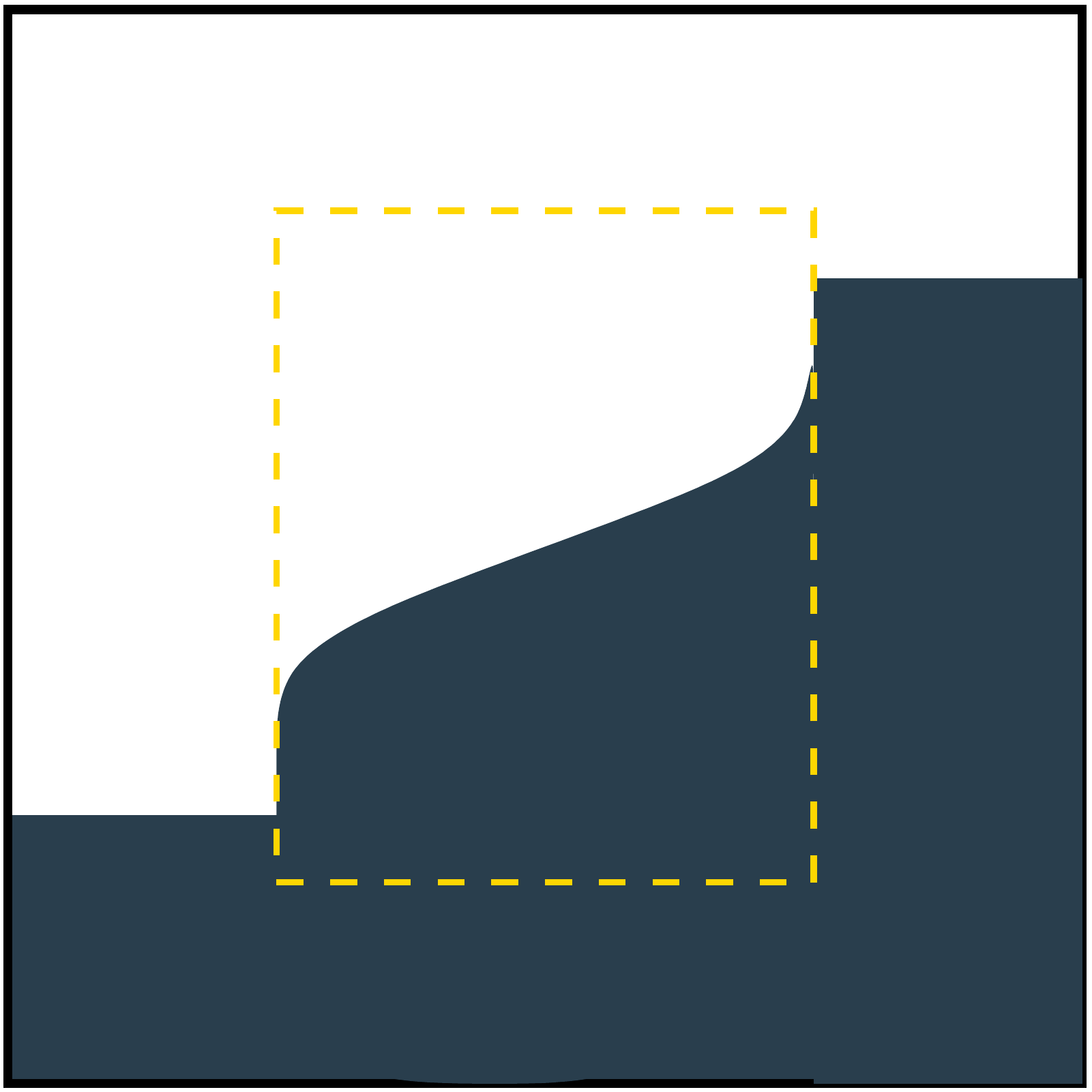}
    \caption{\em {{Classical (on the left) and nonlocal (on the right)
minimal surfaces inside the dashed box, given the outside data.}}}
    \label{STICK}
\end{figure}

Besides its importance in the foundation of a new field of research,
these results contributed in the development of several new lines of investigation:
in particular, the quantity introduced in~\eqref{ZETA} has been later
efficiently utilized in order to detect a rather surprising phenomenon
occurring for the minimizers of the $s$-perimeter, namely their strong tendency
to ``stick'' at the boundary of the domain, in sharp contrast with the classical case, see e.g. Figure~\ref{STICK} -- but this
is somehow a different story.\medskip

Let us now go back to the variational problem
related to the $s$-perimeter. Local minimizers of the $s$-perimeter are called
$s$-minimal sets, and their boundaries are called $s$-minimal surfaces.
An $s$-minimal set which is the subgraph of a function (in a given direction)
is called an $s$-minimal graph. Also,
an $s$-minimal set which is a cone (i.e. a point~$p$ belongs to it
if and only if~$tp$ belongs to it for all~$t>0$) is called an $s$-minimal cone.
The empty set, the full space~$\R^n$, and the half-spaces~$\{\omega\cdot x>0\}$
with~$\omega\in S^{n-1}$
are examples of $s$-minimal cones. The other cones necessarily
exhibit a singularity at the origin, and therefore are called ``singular''. 

In this setting, the regularity of the $s$-minimal surfaces
turns out to be one of the most interesting and challenging topic
related to nonlocal problems, which still presents
many open fundamental questions. 

As usual, Alessio attacked the problem
vigorously, obtaining pioneering results on the topic.
In previous works, in low dimension and for special ranges of the nonlocal exponent,
a regularity in class~$C^{1,\alpha}$ for all~$\alpha\in(0,s)$
was obtained. Then, in~\cite{MR3331523}
a new and general bootstrap result for fully nonlinear
nonlocal equations was provided, whose special
application to the geometric case of $s$-minimal surfaces
established that if an~$s$-minimal surface is locally~$C^{1,\alpha}$ for all~$\alpha\in(0,s)$
(or even just locally Lipschitz, as remarked in~\cite{MR3680376})
then it must be locally~$C^\infty$.

This, combined with previous results, give that
$s$-minimal surfaces are $C^\infty$ in~$\Omega\subset\R^n$ at least in two cases:
\begin{itemize}
\item if~$n\le2$,
\item if~$n\le7$ and~$s\in(s_0,1)$, for some~$s_0\in(0,1)$.
\end{itemize}
It has still to be determined whether or not $s$-minimal surfaces are smooth
in dimension~$n\le7$ for $s$
not close to~$1$ (say, $s$ smaller than the above mention~$s_0$, which,
in principle, could be very close to~$ 1$), and in dimension~$n\ge8$.
In particular, it is still not known an example of singular $s$-minimal cone,
not even in very high dimension.
It is also an open problem to establish whether $s$-minimal surfaces
are analytic.

The regularity of $s$-minimal surface is also related to
the flatness of $s$-minimal graphs. Namely,
we say that the $s$-Bernstein property holds true in~$\R^n$ if
all the $s$-minimal graphs in~$\R^n$
are necessarily affine.
In~\cite{MR3680376}, a general result is given which states that
if there are no singular $s$-minimal cones in dimension~$n-1$,
then the $s$-Bernstein property holds true in~$\R^n$. {F}rom this one obtains that
the $s$-Bernstein property holds true in~$\R^n$
at least in two cases:
\begin{itemize}
\item if~$n\le3$,
\item if~$n\le8$ and~$s\in(s_0,1)$, for some~$s_0\in(0,1)$.
\end{itemize}
Once again, the general picture remains rather mysterious, namely
it is not known whether or not the $s$-Bernstein property holds true in~$\R^n$
when~$n\le8$ and~$s\le s_0$, and when~$n\ge9$, hence we hope that the work of Alessio
will also stimulate new results in these directions.
\medskip

We observe that the $s$-perimeter
defined in~\eqref{INOM} can be also considered in the case~$\Omega=\R^n$,
in which case one simply has that
\begin{equation}\label{INR} {\rm Per}_s(E,\R^n)=I_s(E, E^c).\end{equation}
A well-investigated question in this setting
is the isoperimetric problem, consisting in detecting the minimizers
of the functional in~\eqref{INR} for a prescribed volume. As in the classical
case, these minimizers turn out to be balls, namely, given the scale-invariance
of the $s$-perimeter
\begin{equation}\label{ISOP}
\frac{ {\rm Per}_s(E,\R^n) }{|E|^{\frac{n-s}{n}}}\geq
\frac{ {\rm Per}_s(B_1,\R^n) }{|B_1|^{\frac{n-s}{n}}},
\end{equation}
being~$B_1$ the unit ball of~$\R^n$.

In~\cite{MR3322379}, a number of important questions related to the
$s$-isoperimetric inequality in~\eqref{ISOP}
are addressed. First of all, a ``stable'' version of~\eqref{ISOP}
is obtained, stating that if a set ``almost attains''
the minimal possible value in the fractional isoperimetric
ratio, then it must be necessarily ``almost a ball''.
More precisely, one considers the so-called
Fraenkel asymmetry of a set~$E$,
which measures the~$L^1$-distance of~$E$ from the set of balls of volume~$|E|$
and is defined by
$$ A(E) := \inf_{x\in\R^n}\frac{ |E\Delta B_{r_E}(x)|}{|E|},$$
where~$r_E>0$ is such that~$|B_{r_E}|=|E|$.
In this setting, it is shown in~\cite{MR3322379}
that for any~$s_0\in(0,1)$ and any~$s\in[s_0,1)$ there exists
a positive constant~$C(n,s_0)$ such that
\begin{equation}\label{ISOP2}
\frac{ {\rm Per}_s(E,\R^n) }{|E|^{\frac{n-s}{n}}}\geq
\frac{ {\rm Per}_s(B_1,\R^n) }{|B_1|^{\frac{n-s}{n}}}\left(1+\frac{A^2(E)}{C(n,s_0)}
\right).
\end{equation}
Of course, \eqref{ISOP} is now a particular case of~\eqref{ISOP2}.
Also, the result in~\cite{MR3322379} carries on to the case~$s\nearrow1$.
As usual, the case~$s\searrow0$ is more tricky, and it is conjectured in~\cite{MR3322379} that
$$ C(n,s_0)\sim\frac1{s_0}\qquad{\mbox{ as $s_0\searrow0$}}.$$
The second
variation of the fractional perimeter has been also computed
in~\cite{MR3322379}: this formula can be considered as
a nonlocal counterpart of the classical Jacobi
equation, in which the classical Laplace-Beltrami operator
is replaced by an integral operator along the boundary of the domain,
namely (up to normalization constants) an operator (acting on a given function~$f$)
of the form
$$ \int_{\partial E}\frac{f(x)-f(y)}{|x-y|^{n+s}}\,d{\mathcal{H}}^{n-1}(y),$$
and the norm of the second fundamental form is replaced by
a weighted $L^2$-norm of the normal~$\nu$, such as
$$ \int_{\partial E}\frac{|\nu(x)-\nu(y)|^2}{|x-y|^{n+s}}\,d{\mathcal{H}}^{n-1}(y).$$
Once again,
the nonlocal problems reveal an intrinsic geometric structure
which can be related to the classical objects in the limit as~$s\nearrow1$.
\medskip

Moreover, in~\cite{MR3322379} several variational problems
in which the ``aggregating'' effect of the fractional perimeter
is compensated by a ``disaggregating'' term are
also considered, with special attention to the case given by
Riesz potential~$E$.
In this framework, 
minimizers with small volume are necessarily balls, and the case
of volumes which are not necessarily small demands further investigation.\medskip

Interestingly, minimizers of the $s$-perimeter in~$\R^n$
for a fixed volume satisfy an Euler-Lagrange equation which can be
considered as a prescribed $s$-mean curvature equation. Namely,
for any~$x\in\partial E$, one can consider the $s$-mean curvature
of~$E$ at~$x$, defined by
$$ {\mathcal{H}}^s_E(x):=s(1-s)\int_{\R^n} \frac{
\chi_{E^c}(y)-\chi_{E}(y)
}{|x-y|^{n+s}}\,dy,$$
where, as customary, $\chi_A:\R^n\to\{0,1\}$ denotes the characteristic function of the set~$A$.
By symmetry, one can easily see that the $s$-mean curvature
of a ball of radius~$R$
is constant along its boundary, and, by scaling, it is equal to a constant
depending on~$n$ and~$s$ divided by~$R^s$.

As~$s\nearrow1$, the $s$-mean curvature approaches
the classical mean curvature.

One can show that if~$E$ is a minimizer for the $s$-perimeter in~$\R^n$
for a fixed volume, then its $s$-mean curvature is constant along~$\partial E$.
A natural question in this setting
is to determine the shape of the sets which possess constant $s$-mean curvature
along their boundaries. In the classical setting, this was a classical
result due to the famous Russian mathematician (and
mountaineer)
Aleksandr Danilovich Aleksandrov, stating that 
a smooth, connected and closed hypersurface with constant mean curvature
is necessarily a sphere (hence, soap bubbles are round).

In~\cite{MR3836150} the nonlocal counterpart of the Aleksandrov's
result is obtained, proving that
if a bounded open set with smooth boundary has constant $s$-mean curvature,
then it is necessarily a sphere. It is interesting to remark
that the nonlocal version of such a result is somehow stronger
than the classical case, since the set is not assumed to be
connected. This shows one of the special features of the nonlocal environment,
in which remote interactions give significant contributions to the problem, and, in this
case, they rule out the possibility of multiple connected components (for instance,
two disconnected balls do not have constant $s$-mean curvature).
In other words, the nonlocal setting, in this case, turns out to
be much more rigid than the classical one, since,
even without any connectedness assumption, a set
with constant $s$-mean curvature is a single sphere, whereas of course any
disjoint union of balls with equal radii has constant mean curvature in the classical sense.

Quantitative formulations of 
the Aleksandrov's
result are also provided in~\cite{MR3836150}.
In particular, it is shown that
bounded sets with almost-constant $s$-mean curvature are necessarily
close to a single ball, and moreover
the Lipschitz constant of the $s$-mean curvature
controls the~$C^2$-distance from a single sphere.

Once again, these results reveal some special features
of the nonlocal universe. Indeed, while in the classical case
a connected boundary with almost-constant mean curvature may be close to
a compound of nearby spheres of equal radii, the nonlocal
case turns out to be more rigid
and the quantitative results for almost-constant $s$-mean curvature sets
are obtained without the need of imposing any extra
geometric constraint. 
This also points out an interesting feature of the nonlocal case, 
which prevents bubbling phenomena.

Needless to say, the contributions of Alessio
in this field have been pivotal also to trigger new
research related to sets of constant (and, in particular, zero)
$s$-mean curvature, and to a number of evolution problems
of geometric type (e.g., the ones in which a sets
evolves with normal velocity given by its $s$-mean curvature).\medskip

As customary, Alessio was an avant-garde investigator
of nonlocal geometric and variational problems. His results,
and his research style, will certainly leave an indelible footprint
in the international scenario, and the future research will certainly count on 
his extraordinary talent to solve new questions, open new lines of research
and expand knowledge way beyond the present frontiers.

\end{document}